\documentclass[preprint,11pt]{amsart} 

\usepackage{threeparttable}
\usepackage{caption}
\usepackage{color}
\usepackage{amsmath}
\usepackage{amsfonts}
\usepackage{amssymb}
\usepackage{amsthm}
\usepackage{tikz}

\DeclareMathOperator*{\divergence}{div}
\DeclareMathOperator*{\divergences}{div_s}

\DeclareMathOperator*{\trace}{tr}
\DeclareMathOperator*{\Grad}{\boldsymbol\nabla}
\DeclareMathOperator*{\Grads}{\underline{\boldsymbol{\varepsilon}}}

\DeclareMathOperator{\grads}{\boldsymbol\nabla_s}
\newcommand\RE{\mathbb{R}}
\newcommand\I{\mathbb{I}}
\newcommand{\derivative}[2]{\frac{\partial #1}{\partial #2}}
\newcommand\Huo{H^1_0(\Omega)^d}
\newcommand\Ldo{L^2_0(\Omega)}
\newcommand\Hub{H^1(\B)^d}

\newcommand\Oft{\Omega^f_t}
\newcommand\Ost{\Omega^s_t}
\newcommand\B{\mathcal B}

\newcommand\bfu{\mathbf{u}}
\newcommand\bfv{\mathbf{v}}

\newcommand\n{\mathbf{n}}
\newcommand\w{\mathbf{w}}
\newcommand\x{\mathbf{x}}

\newcommand\s{\mathbf{s}}
\renewcommand\S{\mathbf{S}}
\newcommand\V{\mathbf{V}}
\newcommand\X{\mathbf{X}}
\newcommand\Y{\mathbf{Y}}

\newcommand\LL{\boldsymbol{\Lambda}}
\newcommand\ssigma{\boldsymbol\sigma}
\newcommand\llambda{\boldsymbol\lambda}
\newcommand\mmu{\boldsymbol\mu}

\newcommand\F{\mathbb{F}}
\renewcommand\P{\mathbb{P}}

\newcommand\ds{\mathrm{d}\s}

\newcommand\dt{\Delta t}

\newcommand\T{\mathcal{T}}

\newcommand{\Amatr}{\mathsf{A}}

\newcommand{\Bmatr}{\mathsf{B}}
\newcommand{\Cmatr}{\mathsf{C}}
\newcommand{\Kmatr}{\mathsf{K}}
\newcommand{\Mmatr}{\mathsf{M}}

\newcommand{\Zeromatr}{\mathsf{0}}
\newcommand{\gmatr}{\mathsf{g}}

\newcommand\blue{\begin{color}{black}}
	\newcommand\noblue{\end{color}}
\renewcommand\lg{\begin{color}{black}}
	\newcommand\gl{\end{color}}

\newtheorem{Problem}{Problem}

\begin{document}
	
	\title[A parallel solver for FSI problems with fictitious domain approach]{A parallel solver for FSI problems\\with fictitious domain approach}
	
	\author{Daniele Boffi}
	\address{Computer, Electrical and Mathematical Sciences and Engineering division, King Abdullah University of Science and Technology, Thuwal 23955, Saudi Arabia and Dipartimento di Matematica ``F. Casorati'', Universit\`a degli Studi di Pavia, Via Ferrata 1, 27100, Pavia, Italy}
	\email{daniele.boffi@kaust.edu.sa}
	\urladdr{kaust.edu.sa/en/study/faculty/daniele-boffi}
	\author{Fabio Credali}
	\address{Computer, Electrical and Mathematical Sciences and Engineering division, King Abdullah University of Science and Technology, Thuwal 23955, Saudi Arabia and Dipartimento di Matematica ``F. Casorati'', Universit\`a degli Studi di Pavia, Via Ferrata 1, 27100, Pavia, Italy}
	\email{fabio.credali@kaust.edu.sa}
	\urladdr{cemse.kaust.edu.sa/amcs/people/person/fabio-credali}
	\author{Lucia Gastaldi}
	\address{Dipartimento di Ingegneria Civile, Architettura,Territorio, Ambiente e di Matematica, Universit\`a degli Studi di Brescia, via Branze 43, 25123, Brescia, Italy}
	\email{lucia.gastaldi@unibs.it}
	\urladdr{lucia-gastaldi.unibs.it}
	\author{Simone Scacchi}
	\address{Dipartimento di Matematica ``F. Enriques'', Universit\`a degli Studi di Milano, Via Saldini 50, 20133 Milano, Italy}
	\email{simone.scacchi@unimi.it}
	\urladdr{mat.unimi.it/users/scacchi/}
	
	
	\begin{abstract}
		We present and analyze a parallel solver for the solution of fluid structure interaction problems described by a fictitious domain approach. In particular, the fluid is modeled by the non-stationary incompressible Navier--Stokes equations, while the solid evolution is represented by the elasticity equations. The parallel implementation is based on the PETSc library and the solver has been tested in terms of robustness with respect to mesh refinement and weak scalability by running simulations on a Linux cluster.\medskip
		
	\end{abstract}
	
	\maketitle

\section{Introduction}

\lg The analysis of fluid-structure interaction (FSI) problems is important
for several applications in science and engineering; typical examples that we
have in mind are, for instance, the study of fluid-dynamics of heart valves,
particulate flows, and propagation of free surfaces.
FSI problems are challenging both from the mathematical and computational
point of view: the difficulties originate from the necessity of handling
interactions between several objects and the presence of nonlinear terms in
the governing equations.
As a consequence, during the past years, several approaches have been
presented for the numerical modeling of FSI problems: among those, we mention
the Arbitrary Lagrangian Eulerian formulation (ALE)
\cite{hirt1974arbitrary,donea1977lagrangian,donea2004arbitrary,hughes1981lagrangian},
the unfitted Nitsche method \cite{burman2014unfitted}, the level set
formulation \cite{chang1996level}, and the fictitious domain approach
\cite{glowinski1997lagrange,glowinski2001fictitious}. It is important to
notice that each method is effective for a selected class of problems, since
there is no method that can be applied to all possible situations.\gl

Our fictitious domain approach with distributed Lagrange multiplier was
considered first in \cite{BCG15} as evolution of the immersed boundary method
\cite{peskin2002immersed,pesk}, \lg originally introduced by C. Peskin during
the Seventies for cardiac simulations of heart valves \gl. This approach is
based on the idea that the fluid domain is fictitiously extended also to the
region occupied by the immersed body. In particular, the fluid is governed by
the incompressible time dependent Navier--Stokes equations, while the
structure \lg can be characterized by either linear or nonlinear constitutive
laws for viscous elastic materials. The fluid dynamics is studied on a fixed
Eulerian mesh, while we resort to a Lagrangian description to represent the
motion and the deformation of the immersed structure; the Lagrangian frame is
built by introducing a reference domain which is mapped, at each time step,
into the actual position of the body.\gl

\lg In order to get accurate results, simulations of FSI problems require in
general huge computational resources, both in terms of time and memory. In
this framework, the design of robust parallel solvers is an important tool to
perform, in a reasonable amount of time, computations involving a large number
of time steps and fine space discretizations.
Several works are focused on this task, mainly in the setting of
ALE formulations
\cite{balzani2016,crosetto2011,deparis2016,cai2010,cai2014,heinlein2016,langer2019}
or by applying a variational transfer in order to couple fluid and solids
\cite{krause2016parallel,nestola2019immersed}. To the best of our knowledge,
preconditioners for the fictitious domain formulation with distributed
Lagrange multiplier have been introduced only recently in \cite{boffi2022parallel}.\gl

Here, we continue the analysis of the parallel solver introduced in
\cite{boffi2022parallel}: we focus our attention on the robustness with
respect to mesh refinement, with different choices of time step, and to the weak
scalability. \lg We also describe some implementation issues. \gl
The finite element discretization is performed by choosing the
$(\mathcal{Q}_2,\mathcal{P}_1)$ element for velocities and pressures of the fluid
and the $\mathcal{Q}_1$ element for the structure variables; the time marching
scheme
is a first order semi-implicit finite difference algorithm. Moreover, the
fluid-structure coupling matrix is assembled by exact computations over
non-matching meshes as described in \cite{boffi2022interface}; \lg 
the computational costs of this procedure will be described\gl. At each time
step, the linear system arising from the discretization is solved by
the GMRES method, combined with either a block-diagonal or a block-triangular
preconditioner. \lg Our parallel implementation is based on the PETSc library
from Argonne National Laboratory \cite{petsc-web-page,petsc-user-ref} and on
the parallel direct solver Mumps \cite{amestoy.2001,amestoy.2006}, which is
used to invert the diagonal blocks in the block preconditioners. \gl

%
%
\lg After recalling some functional analysis notation, \gl in
Section~\ref{sec:continuous_formulation} we present the mathematical model
describing fluid structure interaction problems in the spirit of the
fictitious domain approach. In Section~\ref{sec:discrete_formulation}, we
describe the numerical method we implemented for our simulations and in
Section~\ref{sec:parallel_preconditioners} we introduce two possible choices
of preconditioner for our parallel solver. Finally, in
Section~\ref{sec:numerical_tests} we present some numerical tests aiming at
assessing the robustness with respect to mesh refinement and the weak scalability.

\section{Notation}
\lg We recall some useful functional analysis notation \cite{lions2012non}.
Let us consider an open and bounded domain $D$. The space of square integrable
functions is denoted by $L^2(D)$, with scalar product $(\cdot,\cdot)_D$. In
particular, $L^2_0(D)$ is the subspace of functions with null mean over $D$.
We denote Sobolev spaces by $W^{s,q}(D)$: with $s\in\RE$ referring to the
differentiability and $q\in[1,\infty]$ to the summability exponent. When $q=2$,
we adopt the classical notation $H^s(D)=W^{s,2}(D)$. In addition,
$H^1_0(D)\subset H^1(D)$ is the space of functions with zero trace on the
boundary $\partial D$. For vector valued spaces the dimension is explicitly
indicated. \gl

\section{Continuous formulation} \label{sec:continuous_formulation}

We simulate fluid-structure interaction problems characterized by a
visco-elastic incompressible solid body immersed in a viscous incompressible
fluid. We denote by $\Oft$ and $\Ost$ the two regions in $\RE^d$ (with
$d=2,3$) occupied by the fluid and the structure, respectively, at the time
instant $t$; the interface between these two regions is denoted by $\Gamma_t$.
The evolution of such a system takes place inside $\Omega$, that is the union
of $\Oft$ and $\Ost$: this new domain is independent of time and we assume
that it is connected and bounded with Lipschitz continuous boundary
$\partial\Omega$. It is worth mentioning that, even if we are going to consider only thick solids, also the evolution of thin structures can be treated by our mathematical model.

The dynamics of the fluid is studied by considering an Eulerian description,
associated with the variable $\x$. On the other hand, the evolution of the
immersed body is modeled by a Lagrangian description: we introduce the solid
reference domain $\B$, associated with the variable $\s$, so that the
deformation can be represented by the map $\X:\B\longrightarrow\Ost$. This
means that $\x\in\Ost$ is the image \lg at the time $t$ of a certain point
$\s\in\B$ \gl and the motion of the structure is represented by the kinematic
equation
\begin{equation}\label{eq:kinematic}
	\bfu_s(\x,t) = \derivative{\X}{t}(\s,t)\quad\text{for }\x=\X(\s,t),
\end{equation}
denoting by $\bfu_s$ the \lg material \gl velocity. The deformation gradient $\grads\X$ is denoted by $\F$ and $J(\s,t)=\det\F(\s,t)$. \lg Since we are assuming that the solid body is incompressible, the determinant $J$ is constant in time.\gl

In our model, we consider a Newtonian fluid with density $\rho_f$ and viscosity $\nu_f>0$, so that the Cauchy stress tensor can be written as
\begin{equation}
	\ssigma_f = -p_f\I + \nu_f\Grads(\bfu_f),
\end{equation}
where $\bfu_f$ denotes the velocity of the fluid and $p_f$ its pressure \lg
and we denote by $\I$ the identity tensor. \gl In particular, the symbol $\Grads(\cdot)$ refers to the symmetric gradient ${\Grads(\bfv) = (\Grad\bfv +
	\Grad\bfv^\top)/2}$. Therefore, the dynamics in $\Oft$ is governed by the incompressible Navier--Stokes equations
\begin{equation}
	\begin{aligned}
		&\rho_f \bigg( \derivative{\bfu_f}{t} + \bfu_f\cdot\Grad\bfu_f \bigg) = \divergence\ssigma_f\\
		&\divergence\bfu_f = 0.
	\end{aligned}
\end{equation}

For the solid, we consider a viscous-hyperelastic material with density
$\rho_s$ and viscosity $\nu_s>0$; for this type of materials, the Cauchy stress
tensor $\ssigma_s$ can be seen as the sum of two contributions: a viscous
part, similar to the one of the fluid
\begin{equation}
	\ssigma_s^v = -p_s\I + \nu_s\Grads(\bfu_s),
\end{equation}
and an elastic part which can be written\lg, moving from Eulerian to Lagrangian setting, \gl in terms of the Piola--Kirchhoff stress tensor $\P$
\begin{equation}
	\P(\F(\s,t)) = J(\s,t)\ssigma_s^e\lg(\x,t)\gl\F(\s,t)^{-\top} \quad\text{for }\x=\X(\s,t).
\end{equation}
In particular, hyperelastic materials are characterized by a positive energy
density $W(\F)$ which is related with $\P$ since $\P(\F) = \partial
W/\partial\F$. \lg Consequently, the elastic potential energy of the solid
body can be expressed as
\begin{equation}
	E(\X(t)) = \int_{\B} W\big(\F(\s,t)\big) \ds.
\end{equation}\gl
Finally, the system is described by the following equations in strong form
\begin{equation}\label{eq:strong_form}
	\begin{aligned}
		&\rho_f \bigg( \derivative{\bfu_f}{t} + \bfu_f\cdot\Grad\bfu_f \bigg) = \divergence\ssigma_f&&\text{in }\Oft\\
		&\divergence\bfu_f = 0&&\text{in }\Oft\\
		& \rho_s\frac{\partial^2\X}{\partial t^2}=\lg\divergences\gl\big(J\ssigma_s^v\F^{-\top}+\P(\F)\big) &&\text{in }\B\\
		&\divergence\bfu_s = 0 &&\text{in }\Ost
	\end{aligned}
\end{equation}
and completed by two transmission conditions \lg to enforce continuity of velocity and stress \gl along the interface $\Gamma_t$
\begin{equation}\label{eq:transmission}
	\begin{aligned}
		& \bfu_f = \bfu_s&&\text{on }\Gamma_t\\
		&\ssigma_f\n_f = - (\ssigma_s^v + J^{-1}\P\F^\top)\n_s&&\text{on }\Gamma_t,
	\end{aligned}
\end{equation}
where $\n_f$ and $\n_s$ denote the outer normals to $\Oft$ and $\Ost$, respectively.
Moreover, we consider the following initial and boundary conditions
\begin{equation}\label{eq:conditions}
	\begin{aligned}
		&\bfu_f(0) = \bfu_{f,0}&&\text{in }\Omega_0^f\\
		&\bfu_s(0) = \bfu_{s,0}&&\text{in }\Omega_0^s\\
		&\X(0) = \X_0&&\text{in }\B\\
		&\bfu_f = 0&&\text{on }\partial\Omega.
	\end{aligned}
\end{equation}

The idea of the fictitious domain approach is to extend the first two equations in \eqref{eq:strong_form} to the whole domain $\Omega$ \lg so that all the involved variables are defined on a domain which is independent of time. \gl Consequently, following \cite{BCG15}, we introduce two new unknowns
\begin{equation}
	\label{eq:fictitious}
	\bfu=\left\{
	\begin{array}{ll}
		\bfu_f&\text{ in } \Oft\\
		\bfu_s&\text{ in } \Ost
	\end{array}
	\right.\qquad
	p=\left\{
	\begin{array}{ll}
		p_f&\text{ in } \Oft\\
		p_s&\text{ in } \Ost.
	\end{array}
	\right.
\end{equation}
In this new setting, \eqref{eq:kinematic} becomes a constraint on $\bfu$, since we have to impose that
\begin{equation}
	\bfu(\X(\s,t),t) = \derivative{\X}{t}(\s,t)\quad\text{for }\s\in\B.
\end{equation}
This condition can be weakly enforced by employing a distributed Lagrange
multiplier. \lg To this end, we set ${\LL=(\Hub)^\prime}$, the dual space of
$\Hub$, and denote by $\langle\cdot,\cdot\rangle$ the duality pairing between $\LL$ and $\Hub$. Notice that, for $\Y\in\Hub$, we have the following property
\begin{equation}
	\langle \mmu,\Y \rangle = 0 \quad\forall\mmu\in\LL\implies\Y=0.
\end{equation} \gl

At this point, following \cite{BCG15,BG17}, the equations in \eqref{eq:strong_form}, endowed with conditions \eqref{eq:transmission} and \eqref{eq:conditions}, can be written in variational form.
\lg
\begin{Problem}
	\label{pb:pbvar}
	For given $\bfu_0 \in\Huo$ and $\X_0\in W^{1, \infty}(\B)$, find $\bfu(t) \in\Huo$, $p(t) \in\Ldo$, $\X(t) \in\Hub$, and $\llambda(t) \in \LL$ such that for almost all $t \in (0, T)$:
	\begin{equation*}
		\label{eq:FSIvarDLM}
		\left.
		\begin{aligned}
			&\rho_f \left(\derivative{}{t}
			\bfu(t),\bfv\right)_\Omega + b\left(\bfu(t), \bfu(t), \bfv\right) +
			a\left(\bfu(t), \bfv\right)  &&\\
			&\hspace{2.9cm} - \left(\divergence \bfv, p(t)
			\right)_\Omega +\langle \llambda(t), \bfv(\X(\cdot, t)) \rangle = 0 & &\forall
			\bfv \in\Huo \\
			&\left( \divergence \bfu(t), q \right)_\Omega = 0  & &\forall q \in\Ldo
			\\
			&(\rho_s-\rho_f)\left(\derivative{^2 \X}{t^2}(t), \Y \right)_\B +\left( \P(\F(t)),
			\nabla_s \Y\right)_\B  -\langle \llambda(t), \Y \rangle = 0  & &\forall \Y
			\in\Hub \\
			&\langle \mmu, \bfu(\X(\cdot, t), t)- \derivative{\X}{t}(t) \rangle = 0 & &\forall \mmu \in \LL \\
			&\bfu(\x,0) = \bfu_0(\x)  & &\mathrm{in }\ \Omega\\
			&\X(\s,0) = \X_0(\s) & &\mathrm{in }\ \B.
		\end{aligned}
		\right.
	\end{equation*}
\end{Problem}
\gl
In particular, 
\begin{equation}
	\begin{aligned}
		&a(\bfu,\bfv) = \nu\big(\Grads(\bfu),\Grads(\bfv)\big)_\Omega\\ 
		&\lg
		b(\bfu,\bfv,\w)=\frac{\rho_f}2\big((\bfu\cdot\Grad\bfv,\w)_\Omega-(\bfu\cdot\Grad\w,\bfv)_\Omega\big). \gl
	\end{aligned}
\end{equation}
\lg Moreover, \gl $\nu$ is the extended viscosity with value $\nu_f$ in $\Oft$ and $\nu_s$ in $\Ost$. \lg For our numerical tests, we are going to consider $\nu_f=\nu_s$ since it is a reasonable assumption for biological models \cite{peskin1989three}. \gl

For our simulations, we consider a simplified version of the problem: we drop
the convective term of the Navier--Stokes equations and we assume that fluid
\lg and \gl solid materials have the same density, i.e $\rho_s=\rho_f$. \lg We
focus on this case since it is interesting to see how the solver behaves when
this assumption is combined with a semi-implicit time advancing scheme in the
setting of the fictitious domain approach. This is actually the critical
situation when the added mass effect can cause instabilities. For instance, in
\cite{CAUSIN20054506}
a simplified one dimensional setting is considered for which non implicit
schemes are proved to be unconditionally unstable when applied to
FSI problems modeled by ALE
if fluid and solid have the same densities $\rho_s=\rho_f$.
\gl \blue This phenomenon appears regardless of the discrete parameters. In
order to alleviate such critical behavior, some appropriate treatment of the
transmission conditions has been investigated, for instance, in
\cite{deparis2003acceleration,badia2008fluid}) \noblue. \lg
Therefore, the problem we are going to simulate reads as
follows.\gl

\begin{Problem}
	\label{pb:pbvar2}
	Given $\bfu_0 \in\Huo$ and $\X_0\in W^{1, \infty}(\B)$, find $\bfu(t) \in\Huo$, $p(t) \in\Ldo$, $\X(t) \in\Hub$, and $\llambda(t) \in \LL$ such that for almost all $t \in (0, T)$:
	\begin{equation*}
		\label{eq:FSIvarDLM_2}
		\begin{aligned}
			&\rho_f \left(\derivative{}{t}
			\bfu(t),\bfv\right)_\Omega + a\left(\bfu(t), \bfv\right)  &&\\
			&\hspace{1.8cm}- \left(\divergence \bfv, p(t)
			\right)_\Omega +\langle \llambda(t), \bfv(\X(\cdot, t)) \rangle = 0 & &\forall
			\bfv \in\Huo \\
			&\left( \divergence \bfu(t), q \right)_\Omega = 0  & &\forall q \in\Ldo
			\\
			&\left( \P(\F(t)),
			\nabla_s \Y\right)_\B  -\langle \llambda(t), \Y \rangle = 0  & &\forall \Y
			\in\Hub \\
			&\langle \mmu, \bfu(\X(\cdot, t), t)- \derivative{\X}{t}(t) \rangle = 0 & &\forall \mmu \in \LL \\
			&\bfu(\x,0) = \bfu_0(\x)  & &\mathrm{in }\ \Omega\\
			&\X(\s,0) = \X_0(\s) & &\mathrm{in }\ \B.\\
		\end{aligned}
	\end{equation*}
\end{Problem}


\section{Discrete formulation} \label{sec:discrete_formulation}
Before discussing the discrete formulation, we remark that, from now on, we focus on two dimensional problems ($d=2$).

The time semi-discretization of Problem~\ref{pb:pbvar2} is based on the Backward Euler scheme. The time interval $[0,T]$ is \lg uniformly \gl partitioned into $N$ parts with size $\dt=T/N$. We denote the subdivision nodes by $t_n=n\dt$. For a generic function $g$ depending on time, setting $g^n = g(t_n)$, the time derivative is approximated as
\begin{equation}
	\derivative{g}{t}(t_{n+1}) \approx \frac{g^{n+1}-g^n}{\dt}.
\end{equation}
\lg Moreover, the nonlinear coupling terms $\langle \llambda(t), \bfv(\X(\cdot,
t)) \rangle$ and $\langle \mmu, \bfu(\X(\cdot, t), t) \rangle$ are
semi-implicitly treated by considering the position of the structure at the
previous time step as $\langle \llambda^{n+1}, \bfv(\X^n) \rangle$ and $\langle
\mmu, \bfu^{n+1}(\X^n) \rangle$.\gl

For the discretization in space, we work with quadrilateral meshes for both fluid and solid.
For the fluid, we consider a partition $\T_h^\Omega$ of $\Omega$ with meshsize
$h_\Omega$ \lg and \gl two finite element spaces $\V_h\subset\Huo$ and
$Q_h\subset\Ldo$ for velocity and pressure, respectively, satisfying the
inf-sup condition for the Stokes problem. In particular, we work with the
$(\mathcal{Q}_2,\mathcal{P}_1)$ pair, \lg which is one of the most popular
Stokes elements, making use of continuous piecewise quadratic velocities and
discontinuous piecewise linear pressures\gl.

For \lg the solid domain, \gl we choose a partition $\T_h^\B$ of $\B$ with meshsize
$h_\B$, independent of $\T_h^\Omega$. We then consider two finite
dimensional spaces $\S_h\subset\Hub$ and $\LL_h\subset\LL$. We assume that
$\S_h=\LL_h$ and we approximate both \lg the mapping \gl $\X$ and the Lagrange
multiplier \lg $\llambda$ \gl with piecewise bilinear elements on
quadrilaterals. \lg Other stable combinations of finite element spaces for our
class of problems have been studied in \cite{alshehri2022unfitted}, both from
the theoretical and the numerical point of view.\gl

We notice that, since $\LL_h$ is included in $L^2(\B)^d$, at discrete level the \lg duality pairing \gl can be replaced by the scalar product in $L^2(\B)^d$
\begin{equation}
	\lg \langle \mmu_h,\Y_h \rangle \gl = (\mmu_h,\Y_h)_\B \qquad \forall\mmu_h\in\LL_h,\,\forall\Y_h\in\S_h.
\end{equation}

Therefore, we get the following fully discrete problem.
\begin{Problem}\label{pro:full_discretize}
	Given $\bfu_{0,h}\in\V_h$ and $\X_{0,h}\in\S_h$, for all $n=1,\dots,
	N$ find $\bfu_h^n\in\V_h$, $p_h^n\in Q_h$, $\X_h^n\in\S_h$, and $\llambda_h^n\in\LL_h$ fulfilling:
	
	\begin{equation*}
		\label{eq:FSIdiscBDF1}
		\begin{aligned}
			&\rho_f \left( \frac{\bfu_h^{n+1} -
				\bfu_h^{n}}{\dt},\bfv_h\right)_\Omega + a\left(\bfu_h^{n+1}, \bfv_h\right)\\
			&\hspace*{2.2cm} - \left( \divergence \bfv_h, p_h^{n+1}\right)_\Omega+
			\left(\llambda_h^{n+1}, \bfv_h(\X_h^{n})\right)_\B = 0
			& & \forall \bfv_h \in\V_h \\
			&\left( \divergence \bfu_h^{n+1}, q_h \right)_\Omega = 0 
			& & \forall q_h \in Q_h\\
			&\left( \P(\F_h^{n+1}), \nabla_s \Y_h\right)_\B
			- \left(\llambda_h^{n+1}, \Y_h\right)_\B = 0& & \forall \Y_h \in\S_h \\
			&\left(\mmu_h, \bfu_h^{n}(\X_h^{n}) - \frac{\X_h^{n+1} -\X_h^n}{\dt}\right)_\B = 0 & & \forall \mmu_h \in \LL_h \\
			&\bfu_h^0 = \bfu_{0,h}, \quad	\X_h^0 = \X_{0,h}.
		\end{aligned}
	\end{equation*}
\end{Problem}

Assuming for simplicity $\P(\F)=\kappa\F$, Problem~\ref{pro:full_discretize} can be represented in matrix form as
\begin{equation}\renewcommand{\arraystretch}{1.2}
	\label{eq:matrix}
	\left[\begin{array}{@{}cc|cc@{}}
		\Amatr_f		& -\Bmatr^\top	& \Zeromatr			& \Cmatr_f(\X_h^n)^\top \\
		-\Bmatr		& \Zeromatr	& \Zeromatr 			& \Zeromatr \\
		\hline
		\Zeromatr 		& \Zeromatr 	& \Amatr_s   		& -\Cmatr_s^\top \\
		\Cmatr_f(\X_h^n) 	& \Zeromatr   	& -\frac{1}{\dt}\Cmatr_s 	& \Zeromatr \\
	\end{array}\right]
	\begin{bmatrix}
		\bfu_h^{n+1} \\
		p_h^{n+1} \\
		\X_h^{n+1} \\
		\lambda_h^{n+1} \\
	\end{bmatrix}
	= 
	\begin{bmatrix}
		\gmatr_1 \\
		\Zeromatr \\
		\Zeromatr \\
		\gmatr_2 
	\end{bmatrix},
\end{equation}
with
\begin{equation*}
	\begin{aligned}
		& \Amatr_f = \frac{\rho_f}{\dt} \Mmatr_f + \Kmatr_f\\
		& (\Mmatr_f)_{ij} = \left( \boldsymbol{\phi}_j,  \boldsymbol{\phi}_i \right)_\Omega\quad,
		\quad(\Kmatr_f)_{ij} = a\left(\boldsymbol{\phi}_j, \boldsymbol{\phi}_i\right) \\
		& \Bmatr_{ki} = \left( \divergence \boldsymbol{\phi}_i, \psi_k\right)_\Omega\\
		& (\Amatr_s)_{ij} = \kappa \left( \nabla_s \boldsymbol{\chi}_j, \nabla_s \boldsymbol{\chi}_i \right)_\B\\
		& (\Cmatr_f(\X_h^n))_{\ell j} = \left(\boldsymbol{\chi}_\ell, \boldsymbol{\phi}_j(\X_h^n)\right)_\B\quad,
		\quad (\Cmatr_s)_{\ell j} = \left(\boldsymbol{\chi}_\ell, \boldsymbol{\chi}_j\right)_\B\\
		& \gmatr_1 = \frac{\rho_f}{\dt}\Mmatr_f\bfu_h^n,
		\quad \gmatr_2 = -\frac{1}{\dt} \Cmatr_s \X_h^n.
	\end{aligned}
\end{equation*}
Here, $\boldsymbol{\phi}_i$ and $\psi_k$ denote the basis functions of $\V_h$
and $Q_h$ respectively, while $\boldsymbol{\chi}_j$ are the basis functions
of the space defined on $\B$. We observe that, \lg since $\S_h=\LL_h$\gl,
$\Cmatr_s$ \lg is the \gl mass matrix \lg in $\S_h$ or $\LL_h$ \gl. 

We can see that the matrix in \eqref{eq:matrix} splits into four blocks, defined as follows:
\[
\begin{array}{cc}
	\displaystyle \mathcal{A}_{11} = 
	\begin{bmatrix}
		\Amatr_f             & -\Bmatr^\top  \\
		-\Bmatr              & \Zeromatr     
	\end{bmatrix}
	& \displaystyle \mathcal{A}_{12} =
	\begin{bmatrix}
		\Zeromatr                     & \Cmatr_f(\X_h^n)^\top \\
		\Zeromatr                     & \Zeromatr 
	\end{bmatrix}
	\vspace{0.2cm}\\
	\displaystyle \mathcal{A}_{21} =
	\begin{bmatrix}
		\Zeromatr               & \Zeromatr     \\
		\Cmatr_f(\X_h^n)     & \Zeromatr     \\
	\end{bmatrix}
	& \displaystyle \mathcal{A}_{22} =
	\begin{bmatrix}
		\Amatr_s                   & -\Cmatr_s^\top \\
		-\frac{1}{\dt}\Cmatr_s     & \Zeromatr \\
	\end{bmatrix}
	\\
\end{array}
\]
where $\mathcal{A}_{11}$ is related to the fluid dynamic, $\mathcal{A}_{22}$ to the solid evolution, while $\mathcal{A}_{12}$ and $\mathcal{A}_{21}$ contain the coupling term.

Particular attention has to be paid to the assembly of the coupling matrix
$\Cmatr_f(\X_h^n)$, since it involves the integration over $\B$ of solid and
mapped fluid basis functions. \lg In order to compute these integrals, we need
to know how each element $E$ of $\T_h^\B$ is mapped into the fluid domain.
\gl We implement an exact quadrature rule by computing, at each time step, the
intersection between the fluid mesh $\T_h^\Omega$ and a mapped solid \lg
element $\X_h^n(E)$\gl. \lg In order to detect all the intersections, each solid
element is tested against all the fluid elements making use of a bounding box
technique, which in this particular case is trivial since the fluid is discretized with a
Cartesian mesh of squares; then the intersections are explicitly computed by
means of the Sutherland--Hodgman algorithm. \gl For more details about the
procedure in a similar situation, we refer to \cite{boffi2022interface}. \lg
The implementation of this composite rule is quite involved and it is not
straightforwardly parallelizable. \gl

In general, when $\P(\F)$ is nonlinear, we use \lg a solver for nonlinear systems of equations such as the \gl Newton iterator method.

\section{Parallel preconditioners}\label{sec:parallel_preconditioners}
The design of an efficient parallel solver influences two aspects of the
numerical method: first, the finite element matrices need to be assembled in
parallel on each processor, second, the solution of the saddle point
system arising from the discretization has to be solved saving computational
resources, in terms of both memory and execution time. For this purpose, we
implemented a Fortran90 code based on the library PETSc \lg from Argonne
National Laboratory \gl \cite{petsc-web-page,petsc-user-ref}. \lg Such library
is built on the MPI standard and it offers advanced data structures and
routines for the parallel solution of partial differential equations, from
basic vector and matrix operations to more complex linear and nonlinear
equation solvers. In our code, vectors and matrices are built and subassembled
in parallel on each processor. \gl

\lg Our parallel solver adopts \gl two possible choices of preconditioner:
\begin{itemize}
	\item \textit{block-diagonal preconditioner}
	\[
	\begin{bmatrix}
		\mathcal{A}_{11}	& \Zeromatr \\
		\Zeromatr     			& \mathcal{A}_{22} \\
	\end{bmatrix}
	\]
	\item \textit{block-triangular preconditioner}
	\[
	\begin{bmatrix}
		\mathcal{A}_{11}        & \Zeromatr \\
		\mathcal{A}_{21}        & \mathcal{A}_{22} \\
	\end{bmatrix}.
	\]
\end{itemize}

We solve the linear system making use of the parallel GMRES method combined
with the action of our preconditioners, which consists of the exact inversion of the diagonal blocks performed by the parallel direct solver Mumps \cite{amestoy.2001,amestoy.2006}.

\section{Numerical tests} \label{sec:numerical_tests}

The proposed preconditioners have been widely studied in \cite{boffi2022parallel} in terms of robustness with respect to mesh refinement, strong scalability and refinement of the time step. In this work, after reporting new results in terms of optimality, we analyze the weak scalability of our solver. We focus on both linear and nonlinear models describing the solid material.

\lg In the GMRES solver, we adopt as stopping criterion a $10^{-8}$ reduction of the Euclidean norm of the relative residual and a restart parameter of $200$. In the case of the nonlinear model, the stopping criterion adopted for the Newton method is a  $10^{-6}$ reduction of the Euclidean norm of the relative residual. \gl

Our simulations were run on the Shaheen cluster at King Abdullah University of Science and Technology (KAUST, Saudi Arabia). It is a Cray XC40 cluster constituted by 6,174 dual sockets \lg computing \gl nodes, based on 16 core Intel Haswell processors running at 2.3GHz. Each node has 128GB of DDR4 memory running at 2300MHz. 

\subsection{Linear solid model}
We consider a quarter of the elastic annulus ${\{\x\in\RE^2:0.3\le|\x|\le 0.5\}}$ included in $\Omega=[0,1]^2$: in particular, the solid reference domain corresponds to the resting configuration of the body, that is $$\B=\{\s=(s_1,s_2)\in\RE^2:\,s_1,s_2\ge 0,\,0.3\le|\s|\le0.5\}.$$

\begin{figure}
		\begin{center}
			\includegraphics[width=4cm]{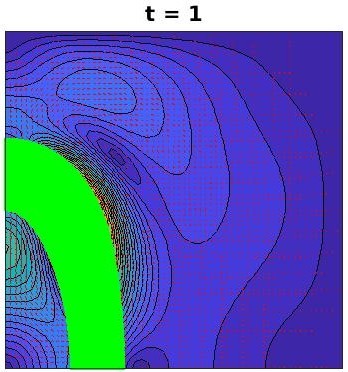}
			\includegraphics[width=4cm]{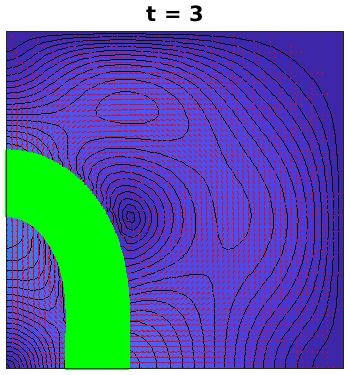}
			\includegraphics[width=4cm]{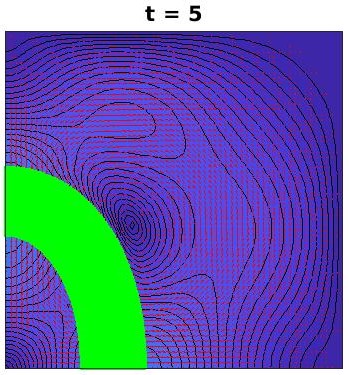}
			\includegraphics[width=4cm]{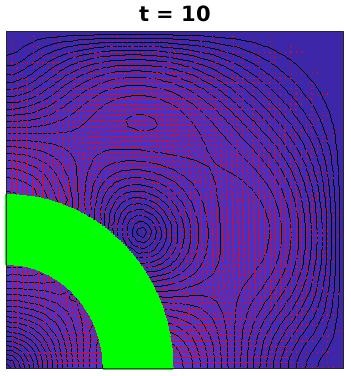}
		\end{center}
		\caption{Four snapshots of the evolution of the structure with linear constitutive law.}
		\label{fig:annulus_evolution}
\end{figure}

The dynamics of the system is generated by stretching the annulus and observing how internal forces bring it back to the resting condition. In this case, $\Omega_0^s$ coincides with the stretched annulus. Four snapshots of the evolution are shown in Figure~\ref{fig:annulus_evolution}.

The solid behavior is governed by a linear model, therefore $\P(\F)=\kappa\,\F$, with $\kappa=10$. We choose fluid and solid materials with same density $\rho_f=\rho_s=1$ and same viscosity $\nu_f=\nu_s=0.1$. We impose no slip conditions for the velocity on the upper and right edge of $\Omega$, while on the other two edges, we allow the motion of both fluid and structure along the tangential direction. Finally, the following initial conditions are considered
\begin{equation*}
	\bfu(\x,0)=0,\qquad\X(\s,0)=\bigg(\frac{s_1}{1.4},1.4\,s_2\bigg).
\end{equation*}

In Table~\ref{tab:linear_optimality}, we report the results for the optimality
test, where the robustness of the solver is studied by refining the mesh and
keeping fixed the number of processors. In particular, we set the time step to
$\dt=0.01$ and the final time of our simulation to $T=2$. The number of processors used for the simulation is $32$. The time $T_{ass}$ needed to assemble the matrix of the problem increases moderately, while the time $T_{coup}$, needed for the assembly of the coupling matrix by computing the intersection between the involved meshes, exhibits a superlinear growth. In terms of preconditioners, we can see that block-diag is not robust with respect to mesh refinement since the number of GMRES iterations grows from $13$ to $430$; clearly, this phenomenon affects also the time $T_{sol}$ we need to solve the system. On the other hand, block-tri is robust since the number of GMRES iterations remains bounded by $14$ when the mesh is refined. Therefore, $T_{sol}$ presents only a moderate growth and, for $1074054$ dofs, it is $30$ times smaller than the value we get for block-diag preconditioner.

\begin{table}
	\centering
		\begin{center}
			\begin{threeparttable}
		\begin{tabular}{r|r|r|r|r|r|r|r|r}
			\hline
			\multicolumn{9}{c}{{\bf Linear solid model -- Mesh refinement test}} \\
			\hline
			\multicolumn{9}{c}{procs = 32, T = 2, $\dt$ = 0.01} \\
			\hline
			dofs    & $T_{ass}(s)$  & $T_{coup}(s)$ & \multicolumn{3}{c|}{block-diag} 		& \multicolumn{3}{c}{block-tri} \\
			&                &               & its           & $T_{sol}(s)$	& $T_{tot}(s)$	& its           & $T_{sol}(s)$	& $T_{tot}(s)$ \\
			\hline
			30534   & 1.02e-2			&9.98e-2 	& 13	& 1.14e-1  	&42.01 		& 7 		&6.93e-2  	&33.83  \\
			120454  & 2.12e-2  			&1.09   	& 31	& 8.30e-1  	&390.90		& 9			&2.40e-1  	&266.17  \\
			269766  & 9.20e-2			&7.60   	& 97	& 5.41 		&2.55e+3 	& 11		&6.47e-1  	&1.65e+3 \\
			478470  & 1.31e-1 			&25.04  	& 192	& 18.75		&9.07e+3 	& 12		&1.14   	&5.24e+3  \\
			746566  & 1.23e-1 			&85.32  	& 422	& 67.92		&3.07e+4	& 13		&2.17   	&1.75e+4  \\
			1074054 & 1.81e-1 			&196.88  	& 430	& 97.19		&5.90e+4	& 14		&3.21   	&4.00e+4  \\
			\hline
		\end{tabular}
	\vspace*{2mm}
	\caption{Refining the mesh in the linear solid model. The simulations are run on the Shaheen cluster. procs = number of processors; dofs = degrees of freedom; $T_{ass}$ = CPU time to assemble the stiffness and mass matrices; $T_{coup}$ = CPU time to assemble the coupling term; its = GMRES iterations; $T_{sol}$ = CPU time to solve the linear system; $T_{tot}$ = total simulation CPU time. The quantities $T_{coup}$, its and $T_{sol}$ are averaged over the time steps. All CPU times are reported in seconds.}
	\label{tab:linear_optimality}
	\end{threeparttable}
	\end{center}
\end{table}	

The weak scalability of the proposed parallel solver is analyzed in
Table~\ref{tab:linear_weak}. Again, we choose $T=2$ and $\dt=0.01$. We perform
six tests by doubling both the global number of dofs and the number of
processors. Thanks to the resources provided by PETSc, the time $T_{ass}$ to
assemble stiffness and mass matrices is perfectly scalable. On the other hand,
the assembly procedure for the coupling matrix is much more complicated: in
order to detect all the intersections between solid and fluid elements, the algorithm consists of two nested loops. For each solid element (outer loop), we check its position with respect to all the fluid elements (inner loop). In particular, only the outer loop is distributed over all the processors. Consequently, $T_{coup}$ is not scalable since the number of fluid dofs, analyzed in serial, increases at each test. We now discuss the behavior of the two proposed preconditioners. It is evident that block-diag is not scalable since the number of GMRES iteration drastically increases as we increase dofs and procs, clearly affecting $T_{sol}$ and $T_{tot}$. On the other hand, block-tri behaves well: even if it is not perfectly scalable, the number of iterations slightly increases from $8$ to $18$ and $T_{sol}$ ranges from $2.24\cdot10^{-1}\,s$ to $11.43\,s$.

\begin{table}
	\begin{center}
		\begin{threeparttable}
		\begin{tabular}{r|r|r|r|r|r|r|r|r|r}
			\hline
			\multicolumn{10}{c}{{\bf Linear solid model -- Weak scalability test}} \\
			\hline
			\multicolumn{10}{c}{T = 2, $\dt$ = 0.01} \\
			\hline
			procs & dofs    & $T_{ass}(s)$  & $T_{coup}(s)$ & \multicolumn{3}{c|}{block-diag} 		& \multicolumn{3}{c}{block-tri} \\
			& &                &               & its           & $T_{sol}(s)$	& $T_{tot}(s)$	& its           & $T_{sol}(s)$	& $T_{tot}(s)$ \\
			\hline
			4 	& 68070	  & 8.55e-2	& 3.95	& 22  & 6.25e-1 & 933.43	& 8	& 2.24e-1	& 833.44 \\
			8 	& 135870  & 1.00e-1 & 5.23	& 38  & 2.16 	& 1.48e+3	& 9	& 4.41e-1	& 1.13e+3 \\
			16 	& 269766  & 1.01e-1 & 8.77	& 111 & 10.23	& 3.80e+3	& 11& 9.70e-1	& 1.95e+3 \\
			32 	& 539926  & 9.24e-2 & 59.27	& 706 & 108.05	& 2.50e+4	& 18& 2.91		& 1.24e+4 \\
			64 	& 1074054 & 1.90e-1 & 48.00	& 429 & 113.59	& 3.24e+4	& 14& 3.90		& 1.04e+4 \\
			128 & 2152614 & 1.90e-1 & 98.63 & -   & -       & -         & 18& 11.43		& 2.20e+4\\
			\hline
		\end{tabular}
		\vspace*{2mm}
		\caption{Weak scalability for the linear solid model. The simulations are run on the Shaheen cluster. Same format as Table~\ref{tab:linear_optimality}.}
		\label{tab:linear_weak}
	\end{threeparttable}
	\end{center}
\end{table}

\subsection{Nonlinear solid model}
\begin{figure}
		\begin{center}
		\includegraphics[width=4cm]{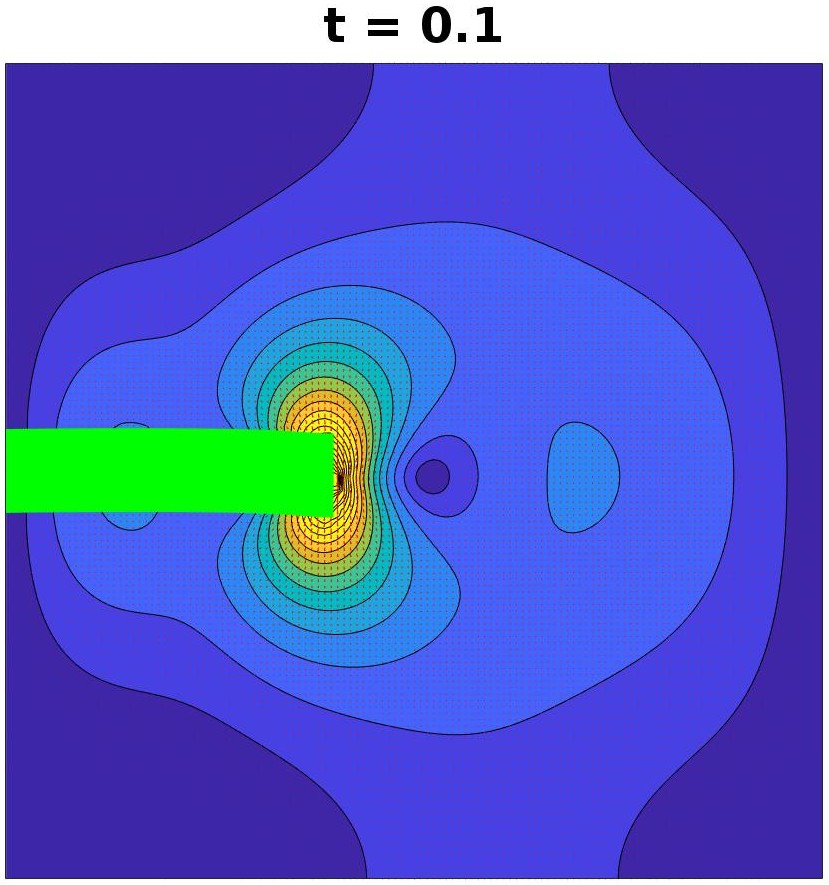}
		\includegraphics[width=4cm]{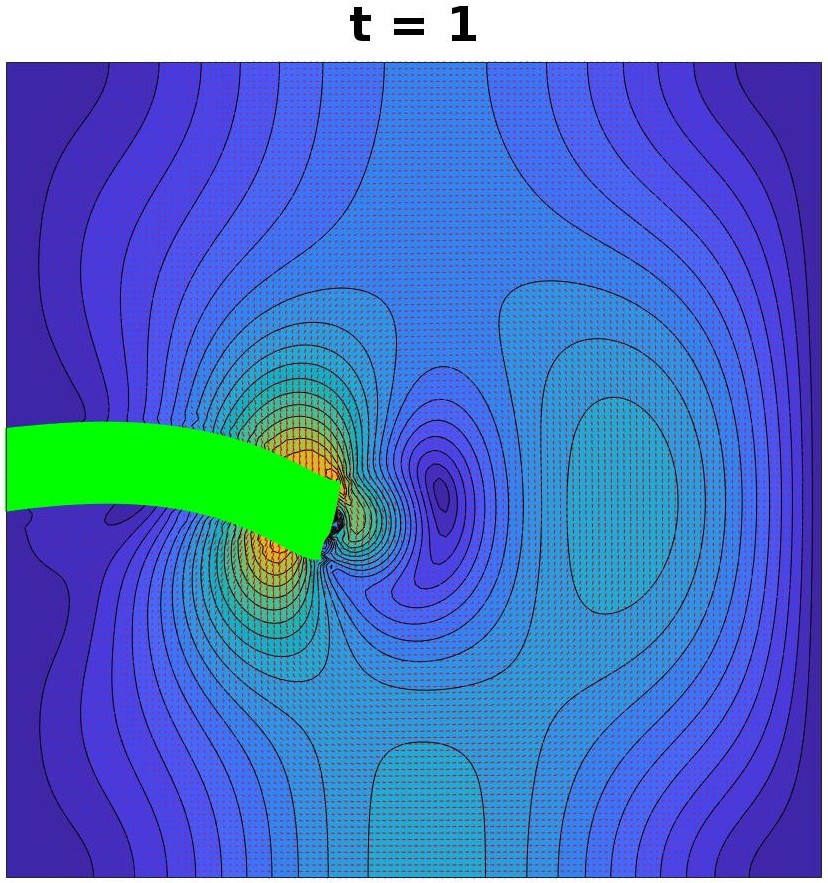}
		\includegraphics[width=4cm]{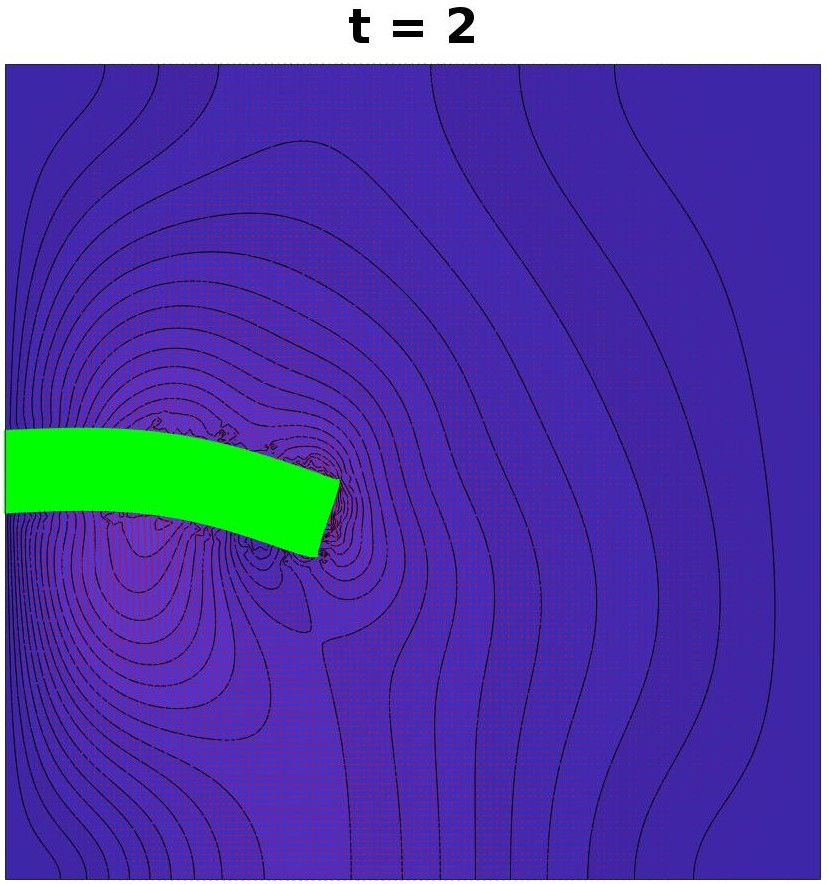}
		\includegraphics[width=4cm]{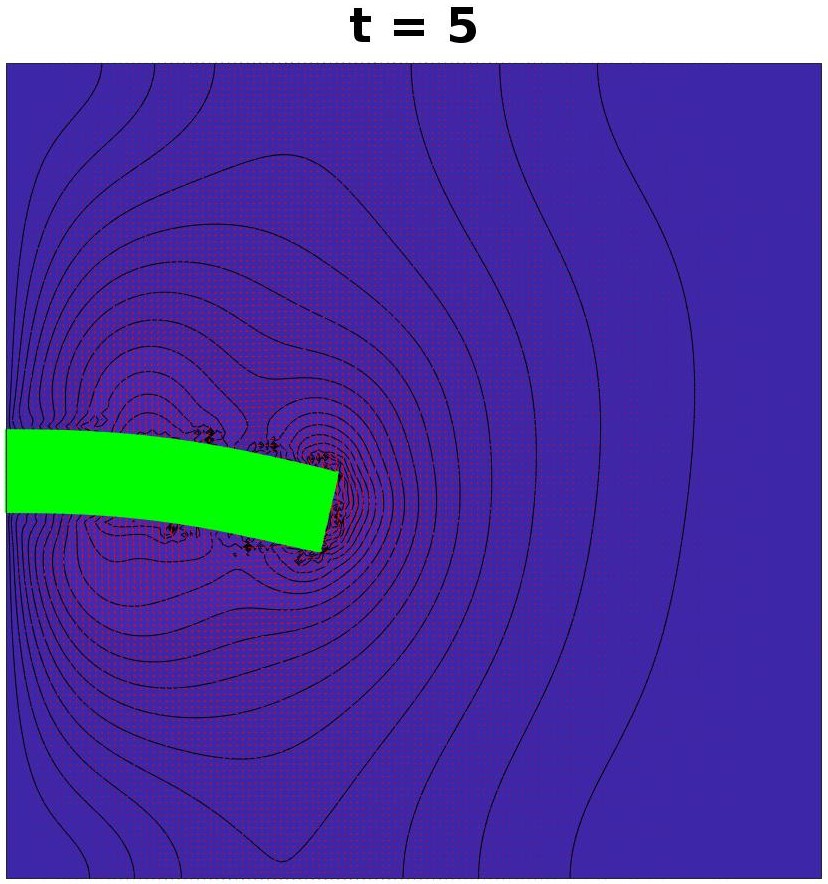}
	\end{center}
	\caption{Four snapshots of the evolution of the structure with nonlinear constitutive law.}
	\label{fig:bar_evolution}
\end{figure}

For this test, we set again the fluid domain $\Omega$ to be the unit square;
the immersed solid body is a bar represented, at resting
configuration, by the rectangle $\B=\Omega_0^s = [0,0.4]\times [0.45,0.55]$.
During the time interval $[0,1]$, the structure is pulled down by a force
applied at the middle point of the right edge. Therefore, when released,
the solid body returns to its resting configuration by the action of internal
forces. Four snapshots of the evolution are shown in
Figure~\ref{fig:bar_evolution}.

The energy density of the solid material is given by the potential strain energy function of an isotropic hyperelastic material; in particular, we have
\begin{equation*}
	W(\F)=(\gamma/2\eta)\exp\big(\eta[\trace(\F^\top\F)-2]\big),
\end{equation*}
where $\trace(\F^\top\F)$ denotes the trace of $\F^\top\F$, while $\gamma=1.333$ and $\eta=9.242$. \lg It can be proved that $W$ is a strictly locally convex strain energy function, as discussed in \cite{delfino1997residual,holzapfel2000new,ogden1997non}.\gl

Also for this test we assume that fluid and solid materials share the same density, equal to $1$, and the same viscosity, equal to $0.2$. The velocity is imposed to be zero at the boundary of $\Omega$, while the following initial conditions are imposed
\begin{equation*}
	\bfu(\x,0)=0,\quad\X(\s,0)=\s.
\end{equation*}

Results for the mesh refinement test are reported in Table~\ref{nonlin_tab_opti_inters}: we consider the evolution of the system during the time interval $[0,2]$, with time step $\dt = 0.002$. The number of processors for the simulations is set to 64, while the number of dofs increases from 21222 to 741702. As for the linear case, $T_{ass}$ increases moderately and $T_{coup}$ follows a superlinear growth. Both preconditioners are robust with respect to mesh refinement: the number of Newton iterations is 2 for each test and the average number of GMRES iterations per nonlinear iteration is bounded by 15 for block-diag and by 10 for block-tri. This behavior of block-diag is in contrast with the results we obtained for the linear solid model: this is due to the finer time step chosen for this simulation.

\begin{table}
		\begin{center}
			\begin{threeparttable}
		\begin{tabular}{r|r|r|r|r|r|r|r|r|r|r}
			\hline
			\multicolumn{11}{c}{{\bf Nonlinear solid model -- Mesh refinement test}} \\
			\hline
			\multicolumn{11}{c}{procs = 64, T = 2, $\dt$ = 0.002} \\
			\hline
			dofs    & $T_{ass}(s)$  & $T_{coup}(s)$ & \multicolumn{4}{c|}{block-diag}               & \multicolumn{4}{c}{block-tri} \\
			&                &               & nit   & its   & $T_{sol}(s)$  & $T_{tot}(s)$  & nit   & its   & $T_{sol}(s)$  & $T_{tot}(s)$ \\
			\hline
			21222   & 4.04e-3  	 & 3.89e-2 	 & 2 	& 11		& 4.29e-1	& 9.35e+2	& 2		& 8   & 3.93e-1	& 8.64e+2 \\
			83398   & 1.68e-2  	 & 3.60e-1 	 & 2 	& 12		& 1.67		& 4.06e+3	& 2		& 8   & 1.57	& 3.86e+3\\
			186534  & 3.80e-2  	 & 1.57 	 & 2 	& 14		& 4.23		& 1.16e+4	& 2		& 9   & 4.00	& 1.11e+4\\
			330630  & 6.68e-2  	 & 4.77 	 & 2 	& 14		& 7.71		& 2.49e+4	& 2		& 10  & 7.07	& 2.37e+4\\
			515686	& 1.05e-1 	 & 11.40 	 & 2 	& 15		& 13.03		& 4.92e+4	& 2		& 10  & 11.48	& 4.58e+4\\
			741702  & 1.52e-1 	 & 23.23 	 & 2  	& 15		& 18.58		& 8.49e+4	& 2 	& 10  &	16.63 	& 7.98e+4\\
			\hline
		\end{tabular}
	\vspace*{2mm}
	\caption{Refining the mesh in the nonlinear solid model. The simulations are run on the Shaheen cluster. procs = number of processors; dofs = degrees of freedom; $T_{ass}$ = CPU time to assemble the stiffness and mass matrices; $T_{coup}$ = CPU time to assemble the coupling term; nit = Newton iterations; its = GMRES iterations to solve the Jacobian system; $T_{sol}$ = CPU time to solve the Jacobian system; $T_{tot}$ = total simulation CPU time. The quantities $T_{coup}$ and nit are averaged over the time steps, whereas the quantities its and $T_{sol}$ are averaged over the Newton iterations and the time steps. All CPU times are reported in seconds.}
	\label{nonlin_tab_opti_inters}
	\end{threeparttable}
	\end{center}
	
\end{table}

In order to study the weak scalability, we choose $T=0.1$ and $\dt=0.002$. The
results, reported in Table~\ref{tab:nonlin_weak}, are similar to the results
obtained for the linear case. As before, $T_{coup}$ is not scalable due to the
algorithm we implemented for the assembling of the coupling term. Even if it is not perfectly scalable, block-tri performs pretty well since the average number of linear iterations per nonlinear iteration increases only from $15$ to $19$. On the other hand, the good behavior of block-diag registered in Table~\ref{nonlin_tab_opti_inters} is not confirmed: the average number of linear iterations reaches $101$, showing a lack of weak scalability as already seen in Table~\ref{tab:linear_weak}.


\begin{table}
	\begin{center}
		\begin{threeparttable}
		\begin{tabular}{r|r|r|r|r|r|r|r|r|r|r|r}
			\hline
			\multicolumn{12}{c}{{\bf Nonlinear solid model -- Weak scalability test}} \\
			\hline
			\multicolumn{12}{c}{T = 0.1, $\dt$ = 0.002} \\
			\hline
			procs & dofs    & $T_{ass}(s)$  & $T_{coup}(s)$ & \multicolumn{4}{c|}{block-diag}               & \multicolumn{4}{c}{block-tri} \\
			&	&                &               & nit   & its   & $T_{sol}(s)$  & $T_{tot}(s)$  & nit   & its   & $T_{sol}(s)$  & $T_{tot}(s)$ \\
			\hline
			4	& 83398	 & 1.01e-1 & 2.21 & 3 & 23	& 6.76 	& 448.65 	& 3	& 15  & 5.50  & 386.13 \\
			8	& 156910 & 1.59e-1 & 3.77 & 3 & 38	& 15.49	& 963.03 	& 3 & 16  & 8.87  & 627.84 \\
			16	& 330630 & 1.62e-1 & 8.92 & 3 & 49	& 36.58	& 2.28e+3	& 3	& 17  & 17.84 & 1.34e+3 \\
			32	& 741702 & 2.60e-1 & 25.18& 3 & 67	& 123.99& 7.46e+3	& 3	& 18  & 48.85 & 3.70e+3 \\
			64	& 1316614& 2.61e-1 & 69.12& 3 & 101	& 328.18& 1.99e+4	& 3	& 19  & 97.28 & 8.38e+3 \\
			\hline
		\end{tabular}
		\vspace*{2mm}
		\caption{Weak scalability for the nonlinear solid model. The simulations are run on the Shaheen cluster. Same format as Table~\ref{nonlin_tab_opti_inters}.}
		\label{tab:nonlin_weak}
		\end{threeparttable}
	\end{center}
\end{table}

\section{Conclusions} \label{sec:conclusions}
We analyzed two preconditioners, block-diagonal and block-triangular, for saddle point systems originating from the finite element discretization of fluid-structure interaction problems with fictitious domain approach. \lg We have focused only on the case where the fluid and solid domains have the same densities, which, based on previous studies, should be the most challenging to solve\gl. In particular, the analysis has been done by studying the robustness with respect to mesh refinement and weak scalability, applying the parallel solver to both linear and nonlinear problems.

Only block-triangular appears to be robust in terms of mesh refinement for linear and nonlinear problems; on the other hand, block-diagonal works well when the time step is very small. 

Moreover, by studying the weak scalability, we can notice two further
limitations of the proposed method, which will be the subject of future
studies. First, the time to assemble the coupling matrix is not scalable: it
is based on two nested loops, related to solid and fluid elements
respectively; but only the external one is done in parallel over the
processors. \lg In order to improve this procedure, one may subdivide the involved meshes into clusters or, alternatively, one may keep track of the position of the solid body at the previous time instant. \gl Second, since the action of the preconditioners consists of the exact inversion of two matrices, the time for solving the linear system slightly increases when the mesh is refined. \lg Some preliminary results have shown that, in contrast with the fluid block $\mathcal{A}_{11}$, the solid block $\mathcal{A}_{22}$ does not behave well when an inexact inversion is applied: our future studies will be primarily focused to overcome this problem. \gl

\lg For what concerns the modeling approach, several extensions are possible.
First of all, we have neglected the convective term in the Navier--Stokes
equation; second, we have considered an isotropic nonlinear
constitutive law for the structure instead of anisotropic variants and,
finally, we have focused only on two dimensional problems.\gl

\section*{Acknowledgments}
The authors are member of INdAM Research group GNCS. D. Boffi, F. Credali and L. Gastaldi are partially supported by IMATI/CNR. Moreover, D. Boffi and L. Gastaldi are partially supported by PRIN/MIUR.

\bibliographystyle{abbrv}
\bibliography{fsi}

\end{document}